\documentclass[fleqn, a4paper, oneside]{imsart}

\usepackage{amsthm,amsmath,amssymb}
\usepackage{graphicx}
\usepackage{makeidx}
\usepackage[authoryear]{natbib}

\bibpunct{(}{)}{;}{a}{,}{,}
\RequirePackage[OT1]{fontenc}

\startlocaldefs

\newtheorem{Definition}{Definition}
\newtheorem{Proposition}{Proposition}
\newtheorem{Theorem}{Theorem}

\newtheorem{Lemma}{Lemma}
\newenvironment{Proof}{\noindent{\bf Proof:}}{\hfill\rule{2mm}{2mm}\vspace{0.3cm}}
\newenvironment{Remark}{\vspace{0.1cm}\noindent{\bf Remark.}}{\vspace{0.3cm}}

\theoremstyle{remark}

\theoremstyle{definition}

\def\R{\mathbb{R}}
\def\C{\mathbb{C}}

\def\P{\mathbb{P}}

\def\eps{\varepsilon}

\def\GG{\mathcal{G}}
\def\TT{\mathcal{T}}
\def\FF{\mathcal{F}}
\def\supp{\mathrm{supp}}

\endlocaldefs

\begin{document}

\begin{frontmatter}


\title{Multiple testing, uncertainty and realistic pictures}
\runtitle{Realistic pictures and uncertainty.}


\begin{aug}
\author{\snm{Mikhail} \fnms{Langovoy}
\ead[label=e1]{langovoy@eurandom.tue.nl}}

\affiliation{
        EURANDOM,\\
         The Netherlands.}

\address{Mikhail Langovoy, Technische Universiteit Eindhoven, \\
EURANDOM, P.O. Box 513,
\\
5600 MB, Eindhoven, The Netherlands\\
\printead{e1}\\
Phone: (+31) (40) 247 - 8113\\
Fax: (+31) (40) 247 - 8190\\}

\and

\author{\snm{Olaf} \fnms{Wittich} \ead[label=e2]{o.wittich@tue.nl}}\corref{}\thanksref{t2}

\affiliation{
        Technische Universiteit Eindhoven and EURANDOM,\\
         The Netherlands.}

\address{Olaf Wittich, Technische Universiteit Eindhoven and \\
EURANDOM, P.O. Box 513,
\\
5600 MB, Eindhoven, The Netherlands\\
\printead{e2}\\
Phone: (+31) (40) 247 - 2499}

\thankstext{t2}{Corresponding author.}

\runauthor{M. Langovoy and O. Wittich}
\end{aug}

\begin{abstract}
We study statistical detection of grayscale objects in noisy images. The object of interest is of unknown shape and has an unknown intensity, that can be varying over the object and can be negative. No boundary shape constraints are imposed on the object, only a weak bulk condition for the object's interior is required. We propose an algorithm that can be used to detect grayscale objects of unknown shapes in the presence of nonparametric noise of unknown level. Our algorithm is based on a nonparametric multiple testing procedure.




We establish the limit of applicability of our method via an explicit, closed-form, non-asymptotic and nonparametric consistency bound. This bound is valid for a wide class of nonparametric noise distributions. We achieve this by proving an uncertainty principle for percolation on finite lattices. \\
\end{abstract}



\end{frontmatter}

\section{Introduction}\label{Section1}


Object detection and image reconstruction for noisy images are two of the cornerstone problems in image analysis. In this paper, we continue our work on an efficient method for quick detection of objects in noisy images. Our approach uses mathematical percolation theory.

Detection of objects in noisy images is the most basic problem of image analysis. Indeed, when one looks at a noisy image, the first question to ask is whether there is any object at all. This is also a primary question of interest in such diverse fields as, for example, cancer detection (\cite{Cancer_Detection_1}), automated urban analysis (\cite{Road_Detection_IEEE}), detection of cracks in buried pipes (\cite{Sinha200658}), and other possible applications in astronomy, electron microscopy and neurology. Moreover, if there is just a random noise in the picture, it doesn't make sense to run computationally intensive procedures for image reconstruction for this particular picture. Surprisingly, the vast majority of image analysis methods, both in statistics and in engineering, skip this stage and start immediately with image reconstruction.

The crucial difference of our method is that we do not impose any shape or smoothness assumptions on the \emph{boundary} of the object. This permits the detection of nonsmooth, irregular or disconnected objects in noisy images, under very mild assumptions on the object's interior. This is especially suitable, for example, if one has to detect a highly irregular non-convex object in a noisy image. This is usually the case, for example, in the aforementioned fields of automated urban analysis, cancer detection and detection of cracks in materials. Although our detection procedure works for regular images as well, it is precisely the class of irregular images with unknown shape where our method can be very advantageous.


We approached the object detection problem as a hypothesis testing problem within the class of statistical inverse problems in spatial statistics. We were able to extend our approach to the nonparametric case of unknown noise density in \cite{langovoy_davies_wittich} and \cite{langovoy_report_Robust_Detection}. In \cite{langovoy_report_2009-035} and \cite{langovoy_davies_wittich}, this density was not assumed smooth or even continuous. It is even possible that the noise distribution is heavy-tailed, see \cite{langovoy_report_2009-035}, \cite{langovoy_davies_wittich} and \cite{langovoy_report_Robust_Detection}.

In \cite{langovoy_wittich_report_R}, we gave an algorithmic implementation of our nonparametric hypothesis testing procedure. We also provided a program that can be used for statistical experiments in image processing. This program was written in the statistical programming language R.

We have shown that there is a deep connection between the spatial structure chosen for the discretisation of the image, the type of the noise distribution on the image, and statistical properties of object detection. These results seem to be of independent interest for the field of spatial statistics.

In our previous papers, we considered the case of square lattices in \cite{langovoy_report_2009-035} and \cite{Langovoy_Wittich_Square}, triangular lattices in \cite{langovoy_davies_wittich} and \cite{langovoy_report_Robust_Detection} and even the case of general periodic lattices in \cite{langovoy_report_Robust_Detection}. In all those cases, we proved that our detection algorithms have linear complexity in terms of the number of pixels on the screen. These procedures are not only asymptotically consistent, but on top of that they have accuracy that grows exponentially with the "number of pixels" in the object of detection. All of our detection algorithms have a built-in data-driven stopping rule, so there is no need in human assistance to stop the algorithm at an appropriate step.

In view of the above, our method can be considered as an unsupervised learning method, in the language of machine learning. This makes our results valuable for the field of  machine learning as well. Indeed, we do not only propose an unsupervised method, but also prove the method's consistency and even go as far as to prove the rates of convergence.


In our previous papers we assumed that the original image was black-and-white and that the noisy image was grayscale. In the present paper, we consider the general case where the signal intensity is completely unknown. This intensity is only assumed to be bounded, but otherwise can vary from pixel to pixel and can be negative.


We propose a multiple testing procedure for detection of grayscale objects of unknown varying intensity in grayscale pictures corrupted by a nonparametric noise that has an unknown distribution. Instead of using a single fixed threshold, we choose a set of thresholds and perform the maximum cluster test from \cite{langovoy_report_Robust_Detection} for each of those thresholds. We show in this paper that, under mild model assumptions, if there is an object in the picture, then it is possible to choose a set of thresholds such that we will consistently detect this object, whenever the object can be even in principle detected on the basis of sizes of percolation clusters. This is one of the two parts that are necessary to prove consistency of the new test.

To establish this result, we need to find out when a signal is too weak so that it \emph{cannot} be detected by our approach. We achieve this goal by proving the so-called uncertainty relation for percolation on finite lattices. This is the main probabilistic result of the present paper. An important distinction of our uncertainty result is that it can be formulated as an explicit condition on the noise distribution and the lattice size. Results of this type are very rare both in statistical literature and in image analysis research. To the best of or knowledge, explicit uncertainty bounds were proved only for Gaussian errors (for example, in research on wavelets by Donoho and coauthors). Our uncertainty relation is much stronger, because it holds uniformly over a wide class of nonparametric error distributions.

Since the problem of detection of greyscale objects cannot be solved in complete generality, we might also provide a set of necessary conditions on the image that makes the object detection possible. We plan to give a possible set of those conditions, as well as the full proof of the consistency theorem for our multiple testing method,  in our forthcoming paper on the subject.



The paper is organized as follows. Section 2 gives a necessary minimal introduction into the mathematical percolation theory. In Section 3, we review our previous results on detection of black-and-white objects in noisy images. In Section 4, we develop an appropriate model for detection of greyscale objects of unknown varying intensity in greyscale pictures corrupted by nonparametric noise. We prove consistency of the basic building blocks of our multiple testing procedure in Theorem 3. The new uncertainty relation for percolation on finite lattices is proved in Section 5. Theorem \ref{uncertainty} of this section is the main mathematical result of the present paper. A new multiple testing procedure for statistical image analysis is proposed in Section 6. Some important results from percolation theory are reviewed in Section 7 of Appendix; this section also contains the proof of the uncertainty relation. Section 8 in Appendix contains the discussion of bounded detector devices.

\section{Percolation theory}


We start with some basic notions of percolation theory. Let $\GG$ be an infinite graph consisting of {\em sites} $s\in\GG$ and {\em bonds} between sites. The bonds determine the topology of the graph in the following sense: We say that two sites $s,s'\in\GG$ are {\em neighbors} if there is a bond connecting them. We say that a subset $C\subset\GG$ of sites is {\em connected} if for any two sites $s,s'\in C$ there are sites $s_1,...,s_n$ such that $s$ and $s_1$, $s_n$ and $s'$, and $s_k$ and $s_{k+1}$ are neighbors for all $k=1,..., n-1$. Considering {\em site percolation}\index{site percolation} on the graph $\GG$ means that we consider random configurations $\omega \in \lbrace 0,1\rbrace^{\GG}$ where the probabilites are {\em Bernoulli}
$$
\begin{array}{ll} P(\omega(s) = 1) = p, & P(\omega (s) = 0) = 1-p
\end{array}
$$
independently for each $s\in\GG$ where $0\leq p \leq 1$ is a fixed probability. If $\omega (s) =1$, we say that the site $s$ is {\em occupied}.\\

\noindent Then, under mild assumptions on the graph, there is a {\em phase transition} in the qualitative behaviour of cluster sizes. To be precise, there is a {\em critical percolation probability}\index{critical probability} $p_c$ such that for $p<p_c$ there is no infinite connected cluster and for $p>p_c$ there is one.\\

\noindent This statement and the very definition $p_c$ being the location of this phase transition are only valid for infinite graphs. We can not even speak of an infinite connected cluster for finite graphs. However, a qualitative difference of sizes of connected clusters of occupied sites can already be seen for finite graphs, say with $\vert\GG\vert = N$ sites. In a sense that will be made precise below, the sizes of connected clusters are typically of order $\log N$ for small $p$ and of order $N$ for values of $p$ close to one. This will yield a criterion to infer whether $p$ is close to zero or close to one from observed site configurations. Intuitively, for large enough values of $N$ the distinction between the two regimes is quite sharp and located very near to the critical percolation probability of an associated infinite lattice.

\section{Maximum cluster test, consistency and rates of convergence}

The signal in our previous papers \cite{langovoy_report_2009-035}, \cite{langovoy_davies_wittich} and \cite{langovoy_report_Robust_Detection} was assumed to be zero-one which corresponds to images with only black and white pixels. In this paper, we will show that the consistency result can be modified to cover also the detection of grayscale objects of unknown intensity. However, first we have to describe our constructions for the basic case.

Let $\GG$ denote a planar graph. We think of the sites $s\in \GG$ as the pixels of a discretized image and of the graph topology as indicating neighboring pixels. In our aforementioned papers, we considered noisy signals of the form
\begin{equation}\label{model}
Y(s) = \bf{1}_{\GG_0} (s) + \sigma\epsilon (s)
\end{equation}
where $\bf{1}_{\GG_0}$ denotes the indicator function of a subset $\GG_0\subseteq \GG$, the noise is given by independent, identically distributed random variables $\lbrace \epsilon (s),s\in \GG\rbrace$ with $E\epsilon = 0$ and $V\epsilon = 1$, and $\sigma > 0$ is the noise variance. Thus, $\sigma^{-1}$ was a measure for the {\em signal to noise ratio}. We refer to \cite{langovoy_report_2009-035} or \cite{langovoy_report_Robust_Detection} for a more detailed introduction.

\begin{Definition}\label{detprob} (The detection problem)\index{detection problem} For signals of the form (\ref{model}), we consider the {\em detection problem} meaning that we construct a test for the following hypothesis
and alternative:
\begin{itemize}
\item[{$\mathbf H_0$}: ] $\GG_0=\emptyset$, i.e. there is no signal.
\item[{$\mathbf H_1$}: ] $\GG_0\neq\emptyset$, i.e. there is a signal.
\end{itemize}
\end{Definition}

In our previous work, we constructed tests for the detection problem given in Definition \ref{detprob} above and computed explicit upper bounds for the {\em type I and type II errors} under some mild condition on the shape of $\GG_0$, called the {\em bulk condition}. We refer to \cite{langovoy_report_2009-035} and \cite{langovoy_report_Robust_Detection} for proofs.

\noindent The setup is as follows: $\TT^{(N)}\subset \TT$ denotes the finite triangular lattice consisting of the $N^2$ sites $s\in\TT$ and bonds which are contained in the subset
$$
\lbrace z\in\C\,:\, \Re(z)\leq N + \frac{1}{2}, \Im(z) \leq \frac{\sqrt{3}}{2}N\rbrace .
$$
By consistency we mean that the test will deliver the correct decision, if the signal can be detected with an arbitrarily high resolution. To be precise, we think of the signal as a subset $G_0\subset\lbrack 0,1\rbrack^2$ and write
$$
G_0^{(N)}:=\lbrace  (N+1/2)x + iN\sqrt{3}y/2\,:\,(x,y)\in G_0\rbrace\subset \C.
$$
The model from equation (\ref{model}) is now depending on $N$, and given by
\begin{equation}\label{model_n}
Y^{(N)}(s) = \mathbf{1}_{\GG_0^{(N)}}(s) + \sigma\,\eps (s)
\end{equation}
where the sites of the subgraph are given by $\GG_0^{(N)} = \lbrace s\in\TT\,:\,s\in G_0^{(N)}\rbrace$
and the bonds of the subgraph are all bonds in $\TT$ that connect two points in $\GG_0^{(N)}$. \\

\noindent We apply now the threshold in the following way. First, we let $\tau = 1/2$, and then define
$$
Y^{(N)}_{\tau}(s) = \left\lbrace\begin{array}{ll} 1 &,Y^{(N)}(s)> 1/2\\
0 & , Y^{(N)}(s)\leq 1/2 \end{array}\right. .
$$

\noindent We consider the following collection of black pixels
\begin{equation}\label{estgnull}
\hat{\GG}_0^{(N)} := \lbrace s\in\TT^{(N)}\,:\,Y^{(N)}_{\tau}(s) =1\rbrace .
\end{equation}

\noindent As a side remark, note that one can view $\hat{\GG}_0^{(N)}$ as an (inconsistent) pre-estimator of $\GG_0^{(N)}$. Now recall that we want to construct a test on the basis of $\hat{\GG}_0^{(N)}$, for the hypotheses $\mathbf{H}_0^{(N)}: \GG^{(N)}_0 =\emptyset$ against the alternative $\mathbf{H}_1^{(N)}: \GG^{(N)}_0 \neq\emptyset$.

\begin{Definition}\label{maxclustest} {\bf (The Maximum-Cluster Test)}\index{maximum cluster test} Let $\phi (N)$ be a suitably chosen threshold depending on $N$.  Let the test statistic $T$ be the size of the largest connected black cluster $C\subset \hat{\GG}_0^{(N)}$.  We reject $\mathbf{H}_0^{(N)}$ if and only if $T \geq\phi(N)$.
\end{Definition}

\noindent For this test, we have the following consistency result under the assumption that the support of the indicator function satisfies the following very weak type of a shape constraint.


\begin{Definition}{\bf (The Bulk Condition)}\index{bulk condition} We say that the support $\GG_0^{(N)}$ of the signal {\em contains a square of side length} $\rho(N)\leq N$ if there is a site $s\in \GG_0^{(N)}$ such that $s+\TT^{(\rho(N))}\subset\GG_0^{(N)}$.
\end{Definition}

\noindent The following consistency result was proved in \cite{langovoy_report_Robust_Detection}.

\begin{Theorem}\label{consistency} For the maximum cluster test, we have
\begin{enumerate}
\item There is some constant $K_0 > 0$ such that for $\phi (N) = K_0\log N$, we have for the {\em type I error}
$$
\lim_{N\to\infty} \alpha (N) = 0.
$$
\item Let $\phi (N)$ be as above. Let the support $\GG_0^{(N)}$ of the signal contain squares of side length $\rho (N)$. If $\rho (N) \geq K_0\log N$, we have for the {\em type II error}


$$
\lim_{N\to\infty} \beta (N) = 0.
$$
\end{enumerate}
In particular, in the limit of arbitrary large precision of sampling, the test will always produce the right detection result.
\end{Theorem}

The next Theorem strengthens Theorem \ref{consistency} and delivers the actual rates of convergence for both types of testing errors. It is a remarkable fact that both types of errors in our method tend to zero exponentially fast in terms of the size of the object of interest. See \cite{langovoy_davies_wittich} or \cite{langovoy_report_Robust_Detection} for the proof.

\begin{Theorem}\label{rates}
Suppose assumptions of Theorem \ref{consistency} are satisfied. Then there are constants $C_1 > 0$, $C_2 > 0$ such that\smallskip

\begin{enumerate}

  \item The type I error of the maximum cluster test does not exceed
  \[
  \alpha (N) \leq \exp(-C_2 \phi (N))
  \]

  \noindent for all $N > \phi (N)$.\medskip

  \item The type II error of the maximum cluster test does not exceed
  \[
  \beta (N) \leq \exp(-C_1 \rho (N))) \,.
  \]

   \noindent for all $N > \rho (N)$.
\end{enumerate}

\end{Theorem}

\section{Realistic pictures}

Instead of the above idealized model, in the present paper we consider the non-distorted signal of interest to be a bounded function $f\in L^{\infty}(\GG)$, i.e. $f(s),s\in\GG$ is a collection of pixel intensities and there exists a $c>0$ such that $\vert f(s)\vert < c$ for all $s\in\GG$. In the sequel, we will call these functions {\em realistic pictures}.\\

\noindent The underlying model for the noisy signal is now as in the indicator signal case given by
\begin{equation}\label{real_model}
Y(s) = f (s) + \sigma\epsilon (s)
\end{equation}
and we assume the same properties of the noise as before in \cite{langovoy_report_Robust_Detection}. More precisely, we assume the following \\

\noindent{\bf Noise Properties.} For a given graph $\GG$, the noise is given by random variables $\lbrace \eps (s)\,:\,s\in\GG\rbrace$ such that
\begin{enumerate}
\item the variables $\eps (s)$ are independent, identically distributed with $E\eps=0$ and $V\eps = 1$,
\item the noise distribution is {\em symmetric},
\item the distribution of the noise is {\em non-degenerate}\index{non-degenerate noise} with respect to a critical probability $p_c$ meaning that if $F$ denotes the cumulative distribution function of the noise and we define
    $$
    \begin{array}{ll}
    m_c^+ = \inf\lbrace x\in\R\,:\,F(x)\geq 1 - p_c\rbrace , & m_c^- = \sup\lbrace x\in\R\,:\,F(x)\leq 1 - p_c\rbrace
    \end{array}
    $$
    then we have $m_c^+ = m_c^-$ where we denote the common value by $m$, and {\em either}
\begin{equation}\label{nondeg1}
F(m) > \lim_{h\to 0, h>0}F(m - h),
\end{equation}
{\em or}
\begin{equation}\label{nondeg2}
F'(m) > 0.
\end{equation}
\end{enumerate}

Furthermore, we assume a {\em bounded detector device} meaning that only signal intensities $\vert Y\vert \leq r$ can be properly displayed, and we assume that this is actually sufficient, i.e. that $\vert Y\vert <r$ for the incoming signal. This is explained more closely in the appendix.


\noindent The test that has to be performed now reads as $\mathbf{H}_0: f=0$ versus the alternative $\mathbf{H}_1: f\neq 0$ where we assume in an analogous way as before, that $f:\lbrack 0,1\rbrack^2 \to \R$ is a bounded continuous function. Thus, in a a similar fashion as before, we construct tests for different resolutions, namely for the hypotheses $\mathrm{H}_0^{(N)}: f^{(N)} = 0$ against the alternatives $\mathrm{H}_1^{(N)}: f^{(N)} \neq 0$ where the discretized function $f^{(N)}:\TT^{(N)}\to\R$ is given by
\begin{equation}
\begin{array}{ll}
f^{(N)}(s) = f(x,y), & s = \left(N + \frac{1}{2}\right)\,x + i\frac{\sqrt{3}}{2}Ny.
\end{array}
\end{equation}
and the corresponding signal is given by
$$
Y^{(N)}(s) = f^{(N)}(s) + \sigma\,\eps (s).
$$
We now have to slightly modify the test, in particular since we do not have any information about the signal strength. This is the main difference to the situation with the indicator function and also the main reason to introduce a {\em bounded detector device}. By that property (and assuming as explained before that the intensity scale provided by the detector is actually sufficient to properly display the signal, or -- likewise -- if we condition on that event) we have a compact scale of thresholds that has to be explored. \\

\noindent Let now $\tau > 0$ and $\GG_{\tau,+}^{(N)}\subset \TT^{(N)}$ denote the {\em super level set}
$$
\GG^{(N)}_{\tau,+} := \lbrace s\in\TT^{(N)}\,:\, Y^{(N)}(s)\geq \tau \rbrace
$$
and $\GG_{\tau,-}^{(N)}\subset \TT^{(N)}$ denote the {\em sub level set}
$$
\GG^{(N)}_{\tau,-} := \lbrace s\in\TT^{(N)}\,:\, Y^{(N)}(s)\leq -\tau \rbrace .
$$
Assume furthermore, that the bounded detector device under consideration has range $r>0$. As a threshold, we use the same $\phi(N) = K_0\,\log N$ as in Theorem \ref{consistency}.\\

\noindent We attempt to do signal detection using the following test statistics
\begin{equation}\label{test_stat}
\begin{array}{lll}T^{(N)}_+(a)& := &\max\,\lbrace \vert C\vert \,:\, C\subset\GG^{(N)}_{a,+} \,\mathrm{black\,\,cluster}\rbrace,\\& &\\
T^{(N)}_-(a)& := &\max\,\lbrace \vert C\vert \,:\, C\subset\GG^{(N)}_{a,-} \,\mathrm{black\,\,cluster}\rbrace\\
\end{array}
\end{equation}
where $a\in \lbrack 0,r/2\rbrack$. It is immediate that we have the following properties as in the case of indicator functions.

\begin{Lemma} Under the null hypothesis, the probabilities that a pixel is erroneously marked black are
\begin{enumerate}
\item $p_{E} = P(s\in \GG^{(N)}_{a,+}) = P(\eps \geq a/\sigma)<1/2 = p_c$,
\item $p_{E} = P(s\in \GG^{(N)}_{a,-}) = P(\eps \leq -a/\sigma)<1/2 = p_c$
\end{enumerate}
and hence subcritical.
\end{Lemma}

\begin{Lemma} {\rm (i)} Let $Q_+^{(N)}\subset\lbrace f^{(N)}\geq a\rbrace$ be a square. Then we have for all $s\in Q_+^{(N)}$ that
$$
p_B = P(s\in \GG^{(N)}_{a/2,+}) = P(\eps \geq -a/2\sigma)>1/2 = p_c .
$$
{\rm (ii)} Let $Q_-^{(N)}\subset\lbrace f^{(N)}\leq -a\rbrace$ be a square. Then we have for all $s\in Q_-^{(N)}$ that
$$
p_B = P(s\in \GG^{(N)}_{a/2,-}) = P(\eps \leq a/2\sigma)>1/2 = p_c .
$$
\end{Lemma}
By these two lemmas, we see that for the test statistics considered above, we are essentially in the same situation as we were in \cite{langovoy_davies_wittich} and \cite{langovoy_report_Robust_Detection}. Both previous lemmas were valid without change if we would consider the respective models
\begin{eqnarray*}
Y^{(N)}_+ &=& \mathbf{1}_{\lbrace f^{(N)}\geq a\rbrace} +
\frac{\sigma}{a}\eps, \\
Y^{(N)}_- &=& \mathbf{1}_{\lbrace f^{(N)}\leq -a\rbrace} +
\frac{\sigma}{a}\eps \\
\end{eqnarray*}
for suitably chosen indicator functions. That implies, we may draw the following conclusion by applying exactly the same proof as in Theorem \ref{consistency}.

\begin{Theorem} For the test statistics considered above, we have:
\begin{enumerate}
\item There is some constant $K_0 > 0$ such that for $\phi (N) = K_0\log N$, we have under the null hypothesis
$$
\lim_{N\to\infty} P(T^{(N)}_+(a)\geq \phi(N)) = \lim_{N\to\infty} P(T^{(N)}_-(a)\geq \phi(N)) = 0.
$$
For $K_0$, we may use the same choice as in Theorem \ref{consistency}.
\item Let $\phi (N)$ be as above. Let $Q_+^{(N)}\subset\lbrace f^{(N)}\geq a\rbrace$ contain squares of side length $\rho (N)$.  If If $\rho (N) \geq K_0\log N$, we have


\[
\lim_{N\to\infty} P(T^{(N)}_+(a/2)\leq \phi(N)) = 0.
\]
\item Let $\phi (N)$, $\rho(N)$ be as above. Let $Q_-^{(N)}\subset\lbrace f^{(N)}\leq -a\rbrace$ contain squares of side length $\rho (N)$.  Then we also have
\[
\lim_{N\to\infty} P(T^{(N)}_-(a/2)\leq \phi(N)) = 0.
\]
\end{enumerate}
In particular, the test statistic associated to the correct scale parameter $a/2$ will asymptotically always produce the right detection result.
\end{Theorem}

\noindent At first sight, the situation seems rather similar as for indicator functions in Theorem \ref{consistency}. However, it is completely different in the sense that the consistency result only holds if we pick the right signal strength in advance. We might be able to overcome this problem by considering a scale of tests for some positive $a > 0$.

\section{Uncertainty}

It is intuitively clear that, for principal reasons, it is not possible to detect a signal with arbitrarily small signal to noise ratio on a lattice of finite size, no matter which method is used for detection. However, for every particular method, it is very difficult to provide a "horizon of consistency" in explicit form. Results of this type a very rare in hypothesis testing, image analysis or machine learning. Typically, one proves those results in special cases like Gaussian errors. 

In this section, we provide an explicit, closed-form, non-asymptotic and nonparametric consistency bound for our method. This bound is valid for a wide class of nonparametric noise distributions and is given in Theorem \ref{uncertainty}.


Recall from the proof of Theorem \ref{consistency} that the constant $K_0$ in the threshold was given by the inequality
$$
K_0 = 2C > \lambda (p_E)^{-1}
$$
where $p_E$ is the (subcritical) probability under the null hypothesis that a pixel is erroneously marked black and $\lambda (p_E)$ is the constant from the Aizenman - Newman theorem, see \cite{langovoy_report_Robust_Detection}. Thus, we have to begin by finding a proper estimate of $\lambda (p)$ for a subcritical $p$.\\


\noindent The classical Aizenman-Newman theorem reads as follows.\\

\begin{Proposition}\label{AizNew}{\bf (Aizenman-Newman Theorem)}\index{Aizenman - Newman theorem} Consider percolation with subcritical probability $p < p_c = 1/2$ on the infinite triangular lattice $\TT$. Then there is a constant $\lambda (p) > 0$ depending on $p$
such that
\begin{equation}\label{subcrit}
 P(\vert C\vert \geq n) \leq e^{-n\,\lambda(p)}
\end{equation}
for all $n\geq 1$ where $C$ denotes the black cluster containing the origin.
\end{Proposition}

\begin{Remark} Please note that we use asymptotics-oriented estimates to prove statements about finite lattices. For instance in the case of (\ref{EXPECTBOUND}) below, these estimates are not the best possible. So we can by no means expect that the bound in Theorem \ref{uncertainty} is sharp. But it is good enough to serve as an illustration of the basic principle.
\end{Remark}

\noindent In the sequel, $\chi (p)$ denotes the expected size of the cluster containing $0\in\TT$ in the infinite lattice depending on the subcritical occupation probability $p<p_c=1/2$.

\begin{Lemma} For the infinite triangular lattice, we have
\begin{equation}\label{EXPECTBOUND}
\chi (p) \leq \frac{1}{18} \,\vert p - p_c\vert^{-1}.
\end{equation}
\end{Lemma}

\begin{Proof} See Appendix \ref{pf}.
\end{Proof}

\noindent By the Aizenman - Newman Theorem (Proposition \ref{AizNew}), we obtain an upper bound for this expectation value by
\begin{equation*}
\chi (p) = \sum_{n\geq 1} P(\vert C\vert \geq n)\leq \sum_{n\geq 1} e^{-n\lambda (p)} = \frac{e^{-\lambda (p)}}{1 - e^{-\lambda (p)}}.
\end{equation*}
Thus, we have
\begin{equation}\label{Lambda}
\lambda (p) \leq - \log \left( \frac{\chi(p)}{1+\chi(p)}\right) = \log\left(1 + \frac{1}{\chi(p)}\right).
\end{equation}
Combining these two results yields

\begin{Lemma}\label{Lam2} We have
\begin{equation*}
\lambda(p)^{-1} \geq \frac{1}{\log\left(1 + 18\,\vert p - p_c\vert\right)}.
\end{equation*}
\end{Lemma}

\begin{Proof} (\ref{EXPECTBOUND}) together with (\ref{Lambda}) implies
$$
\lambda (p) \leq \log \left( \frac{1+\chi(p)}{\chi(p)}\right) = \log \left( 1 + \frac{1}{\chi(p)}\right)\leq\log\left(1 + 18\,\vert p - p_c\vert\right) .
$$
\end{Proof}

\noindent This implies that an intrinsic feature of the procedure is the following form of {\em uncertainty}: {\em By our procedure, we cannot -- even in principle -- detect signals with arbitrary low signal to noise ratio on a finite lattice of fixed size.} \\

\noindent Of course, it is clear that something like the above statement is valid for any statistical testing procedure. Therefore, something similar is also valid for signal detection. An important distinction of our uncertainty result is that we can give an explicit condition on the noise level and the lattice size, such that this condition implies that our test \emph{does not} work. Results of this type are very rare both in statistical literature and in image analysis research. To the best of or knowledge, explicit uncertainty bounds were proved only for Gaussian errors (for example, in research on wavelets by Donoho and coauthors). Our uncertainty relation is much stronger, because it holds irrespectively of the actual noise distribution, uniformly over a wide class of nonparametric error distributions.

To be precise, we consider again the threshold $\phi (N) = K_0\,\log N$ and the fact that in the proof of Theorem \ref{consistency}, we had to choose $K_0 = 2C > 2\lambda(p_E)^{-1}$. That implies together with Lemma \ref{Lam2} that
\begin{equation}\label{Lem3}
\phi (N) = K_0 \,\log (N) > \lambda(p)^{-1}\,\log N^2 \geq \frac{\log N^2}{\log\left(1 + 18\,\vert p_E - p_c\vert\right)}.
\end{equation}
But that means, that for values of $p$ which are very close to the critical probability, the threshold $\phi (N)$ may exceed the lattice site $N^2$ and our method breaks down. To be precise, we have the following statement.

\begin{Proposition}\label{uncertainty2} If the lattice size $N^2$ is fixed, the threshold $\phi(N)$ is larger than the lattice size, and therefore, the null hypothesis will never be rejected, if we have
$$
\vert p_E - p_c\vert < \frac{1}{18}\, \left\lbrace\left( N^2\right)^{\frac{1}{N^2}} - 1\right\rbrace.
$$
\end{Proposition}

\begin{Proof} By (\ref{Lem3}), we have $\phi(N) > N^2$ if
$$
\frac{\log N^2}{\log\left(1 + 18\,\vert p_E - p_c\vert \right)} > N^2.
$$
\end{Proof}

\noindent Finally, we want to relate this statement directly to the signal strength. Thus, if $\vert f^{(N)} \vert \leq a$ we say that the {\em signal to noise ratio} is given by $\rho = a/\sigma$. Let us now assume that the distribution function of the noise $F$ is continuous at zero. Then
$$
\vert p_E - p_c\vert = \frac{1}{2} - p_E = P(0 < \eps < a/\sigma) = F(\rho) - F(0)
$$
and we finally obtain

\begin{Theorem}\label{uncertainty} {\bf (Uncertainty)}\index{uncertainty} Assume that the distribution function of the noise is continuous at zero. A signal $f^{(N)}$ with $\vert f^{(N)} \vert \leq a$ and signal to noise ratio $\rho = a/\sigma$ can only be detected by our method, if
\begin{equation*}
\frac{P(0 < \eps < \rho)}{\left( N^2\right)^{\frac{1}{N^2}} - 1} > \frac{1}{18}
\end{equation*}
that means if either the lattice size is sufficiently large or the signal to noise ratio is sufficiently small.
\end{Theorem}

\begin{Remark} {\rm (i)} Note that this statement only means that -- as a matter of principle -- we can not  detect signals of arbitrarily small strength on a finite lattice of a given size. That does not at all mean that detection of signals that respect the bound above is automatically possible in an effective way. Topics like type I and type II error are not at all touched by the uncertainty bound. In other words, from the uncertainty relation we can derive only a necessary condition for the signal to be detectable via our method. Usually this condition will not be sufficient.\\

\noindent {\rm (ii)} From studying the behavior of the function
$$
s(x) = \frac{1}{18}\left( e^{-x\ln x}-1\right)
$$
on the unit interval, we see that the bound is always fulfilled if $P(0<\eps < \rho) > 0.25 \approx \max_{x\in\lbrack 0,1\rbrack} s(x)$.
\end{Remark}

\noindent The proof of Theorem \ref{uncertainty} consists of a simple reformulation of the preceding proposition and is therefore omitted. However, we still have to justify why we use the word {\em uncertainty} in connection with this statement. This discussion can only be purely informal. The analogy simply is that a function of the signal to noise ratio times another function of the lattice size have to exceed a certain value for a signal to be detectable. Otherwise, the signal is virtually not existing. A weaker version of the statement might provide another argument: There is some number $M > 0$ such that
\[
s(x) \leq M\sqrt{x}
\]
for all $x\in\lbrack 0,1\rbrack$. If we assume now that $F$ has a continuous and sufficiently smooth density $f$ with $f''(0) <0$, we have the weaker statement that the signal can be detected only if
\[
f(0)\rho \geq F(\rho) - F(0) > M\sqrt{1/ N^2}
\]
or, if
\begin{equation}\label{simple_uncertainty}
N\,\rho > \frac{M}{f(0)}\,.
\end{equation}
Thus, for a signal to be detectable, the product of two conjugate parameters may not exceed a bound given by the circumstances. Otherwise, the signal is not detectable, even in principle.\\

\section{Multiple testing for realistic pictures}

By the uncertainty principle, we obtain a minimal threshhold value below which it does not make any sense to try to detect a signal. So there is a natural lower bound $\tau_0$ for a threshold. The upper bound is provided by the size $r$ of the bounded detector device. That means, the range of intensities of detectable signals is $\lbrack -r,-\tau_0\rbrack \cup \lbrack \tau_0,r\rbrack$. Thus, if $f$ is the signal, and we assume bulk conditions for the super-level sets as in Corollary 1, taking into account the the simple {\em monotonicity property} that $a > a'$ implies ${\mathbf{1}}_{\lbrace f^{(N)}\geq a\rbrace}\leq{\mathbf{1}}_{\lbrace f^{(N)}\geq a'\rbrace}$, we will certainly be able to consistently detect an object (if the object can be potentially detected on the basis of percolation clusters), via the following scheme:
\begin{enumerate}
 \item Consider the threshold scheme
$$
a_k =   2^{-k} r , k=1,...,N.
$$
\item Beginning with $a=a_{1}$, calculate the test statistics $T_-^{(N)}(a)$, $T_+^{(N)}(a)$. Terminate, if either the null hypothesis is rejected (at a properly adjusted level, if necessary) or if you reach $a_k$ with $k \geq \log (r/\tau_0)$.
\end{enumerate}

It can be shown that, under certain conditions on $f$ and $\sigma$, we would have to repeat the maximum cluster test at most $O (\log N)$ times. Since each repetition of the maximum cluster test takes $O(N^2)$ operations, the new multiple testing procedure is going to take at most $O (N^2 \log N)$ operations overall. Since the input size is $N^2$, this implies that under some conditions our initial procedure (of linear complexity) slows down by a logarithmic factor. Asymptotically, this is only slightly slower than the original test, but the new test is adaptive with respect to the unknown image color intensity.

A point that needs to be addressed carefully here is the probability of false rejection of the null hypothesis. Indeed, we perform here not a single test, but a collection of up to $O (\log N)$ tests, and results of those tests are \emph{not} independent. This is a basic question that always occurs in the field of multiple testing. Luckily for us, for each of the thresholds $a_k$ the direct analog of Theorem \ref{rates} holds: the type I error of any single test tends to zero exponentially fast, while the power tends exponentially fast to one. Moreover, our tests are "monotonous" with respect to the threshold value (since the maximum cluster size is an increasing event). This also implies that we have to pay attention only to those thresholds $a_k$ where at least one of the level sets $\GG^{(N)}_{\tau,-}$ and $\GG^{(N)}_{\tau,+}$ doesn't contain black clusters crossing the whole screen. Using those properties, we will be able to combine the results of not more than $O (\log N)$ tests $T_-^{(N)}(a_k)$ and $T_+^{(N)}(a_k)$ and get a unique decision out of them, while keeping the type I error of the multiple test controlled. We plan to present those results in succeeding papers.

\smallskip
\noindent {\bf Acknowledgments.} The authors would like to thank the EURANDOM Report Series reviewers for carefully reading this manuscript. \\

\bibliographystyle{plainnat}

\bibliography{papiere}

\bigskip
\noindent {\bf Appendix.}\\

\section[Some facts from percolation theory]{Some facts from percolation theory}

\noindent In this section, we collect some basic statements and techniques from the theory of percolation. In particular, we are going to prove the inequality (\ref{EXPECTBOUND}) which is basic for the introduction of {\em uncertainty principle}.\\

\subsection{FKG and BK inequality}

Recall the partial ordering
$$
\begin{array}{ll} \omega_1 \preceq \omega_2 & :\Longleftrightarrow \,\, \omega_1(s)\leq\omega_2(s)\,\,\mathrm{for \,\, all}\,\, s\in\TT .   \end{array}
$$
on the set $\Omega = \lbrace 0,1\rbrace^{\TT}$ of all {\em percolation configurations} from Definition \ref{increase} and that an {\em event} $A\subset \Omega$ is {\em increasing} if we have
$$
\mathbf{1}_A (\omega_1) \leq \mathbf{1}_A (\omega_2)
$$
for the corresponding indicator variable whenever $\omega_1\preceq\omega_2$. \\

\noindent The FKG inequality was already stated before and is just added here another time for completeness.

\begin{Proposition}{\bf (FKG inequality)} If $A$ and $B$ are both increasing (or both decreasing) events, then we have
$$
P(A\cap B) \geq P(A)\, P(B).
$$
\end{Proposition}

\begin{Proof}\cite{FoKaGi:71}\end{Proof}

\noindent Let $\GG\subset \TT$ be a {\em finite} subgraph and
$$
\FF_{\GG}\sigma (\lbrace 0,1\rbrace^{\GG})\subset \sigma (\lbrace 0,1\rbrace^{\TT}) =:\FF_{\TT}
$$
the sub sigma - algebra associated to the percolation configurations on $\GG$ (in the canonical version).
Let now $A,B\in\FF_{\GG}$ be two increasing events. We define the support of $\omega\in \lbrace 0,1\rbrace^{\GG}$ to be $$
\mathrm{supp} \omega := \lbrace s\in\GG\,:\,\omega (s) = 1 \rbrace .
$$
and for a subset $H\subset \supp\omega$, we write
$$
\omega\vert_H := \left\lbrace\begin{array}{ll} 1 & s\in H \\ 0 & {\mathrm else}\end{array}\right. .
$$

\begin{Definition} Let $A$, $B$ be as above. The event $A\circ B$ that $A$ {\em and} $B$ {\em occur disjointly} is given by
$$
A\circ B := \lbrace\omega\in\lbrace 0,1\rbrace^{\TT}\,:\,\exists_{H(\omega)\in\supp \omega}\,\omega\vert_{H(\omega)}\in A, \omega\vert_{\supp \omega - H(\omega)}\in B \rbrace .
$$
\end{Definition}

\noindent The BK inequality now reads as follows.

\begin{Proposition}{\bf (BK inequality)} Let $A,B\in\FF_{\GG}$ be increasing events. Then
$$
P(A\circ B) \leq P(A)\,P(B) .
$$
\end{Proposition}

\begin{Proof}\cite{Gri:99}, p. 38 ff.\end{Proof}

\subsection{Russo's formula}

\noindent Let $s\in\TT$ be a site. We consider the involution $j_s: \lbrace 0,1\rbrace^{\TT}\to\lbrace 0,1\rbrace^{\TT}$ given by
$$
j_s (\omega) (s') := \left\lbrace\begin{array}{ll} \omega (s') & s'\neq s\\ 1 - \omega (s') & s' = s\end{array}\right. .
$$
From this definition, we see that the configuration space is a disjoint union
$\lbrace 0,1\rbrace^{\TT} = \Omega(s)^+ \cup j_s\Omega(s)^+$, where
$$
\Omega(s)^+ := \lbrace \omega\in\lbrace 0,1\rbrace^{\TT}\,:\, \omega(s) = 1\rbrace.
$$

\begin{Definition} {\bf (Pivotal sites)}\index{pivotal sites} Let $\GG\subset\TT$ be a finite subgraph and $A\subset\FF_{\GG}$ be an increasing event. The event {\em the site} $s$ {\em is pivotal for} $A$ is given by
$$
\mathrm{Piv}(A,s) := \lbrace \omega\in\lbrace 0,1\rbrace^{\TT}\,:\,\mathbf{1}_A(\omega) \neq \mathbf{1}_A\circ j_s(\omega)     \rbrace .
$$
\end{Definition}

\noindent Russo's formula is a statement about how the probability of a certain event changes if the individual site occupation probability $p$ is changed. We denote by $P_p(A)$ the probability of the event $A$ if this probability is $p$ and by
$$
N_A := \sum_{s\in\GG} \mathbf{1}_{\mathrm{Piv}(A,s)}
$$
the number of pivotal elements for $A$.

\begin{Proposition}{\bf (Russo's formula)}\index{Russo's formula} Let $\GG\subset\TT$ be a finite subgraph and $A\subset\FF_{\GG}$ be an increasing event. Then
\begin{equation}\label{russo}
\frac{d}{dp} P_p(A) = E_pN_A .
\end{equation}
\end{Proposition}

\begin{Proof} {\rm (i)} First of all, since $A$ is increasing and $\omega \preceq j_s(\omega)$ for all $\omega\in\Omega(s)^+$, we have on the set $\Omega(s)^+$
$$
\left(\mathbf{1}_A - \mathbf{1}_A\circ j_s\right)(\omega) = \left\lbrace\begin{array}{ll}
1 & \omega\in\mathrm{Piv}(A,s)\cap\Omega(s)^+\\
0 & \mathrm{else}\end{array}\right. =  \mathbf{1}_{\mathrm{Piv}(A,s)\cap\Omega(s)^+} .
$$
{\rm (ii)} In the sequel, we write $\Omega(s)^- = j_s\Omega(s)^+$. Let $\mathbb{P}\vert_{\Omega(s)^-}$ denote the restriction of the measure to $\Omega(s)^-$. Then, the image measure under $j_s$ is a measure on $\Omega(s)^+$ with density
$$
\frac{d\mathbb{P}\vert_{\Omega(s)^-}\circ j_s}{d\mathbb{P}} = \frac{\mathbb{P}(\Omega(s)^+)}{\mathbb{P}(\Omega(s)^-)}.
$$
That implies
\begin{eqnarray*}
E_p\left\lbrack\mathbf{1}_A\vert\Omega(s)^{-}\right\rbrack &=& \frac{\int_{\Omega(s)^{-}} \mathbf{1}_A (\omega)\mathbb{P}(d\omega)}{\mathbb{P}(\Omega(s)^-} =  \frac{\int_{\Omega(s)^{+}}\mathbf{1}_A\circ j_s (\omega')\mathbb{P}\circ j_s(d\omega')}{\mathbb{P}(\Omega(s)^-)} \\  &=& E_p\left\lbrack\mathbf{1}_A\circ j_s\vert\Omega(s)^{+}\right\rbrack .
\end{eqnarray*}
{\rm (iii)} Now let $p'>p$ and denote by $E_{p',s}$ the expectation with respect to the product measure $\P_s$ with marginals
$$
\P_s (\omega (s') = 1) = \left\lbrace \begin{array}{ll} p' & s' = s \\ p & \mathrm{else} \end{array}\right. .
$$
Thus
\begin{eqnarray*}
& & P_{p',s}(A) - P_p(A) = E_{p',s}\mathbf{1}_A - E_p\mathbf{1}_A \\
&=& P_{p',s}(\Omega(s)^+)E_{p',s}\left\lbrack\mathbf{1}_A\vert \Omega(s)^+\right\rbrack  + P_{p',s}( \Omega(s)^-)E_{p',s}\left\lbrack\mathbf{1}_A\vert \Omega(s)^-\right\rbrack \\
& & - P_{p}(\Omega(s)^+)E_{p}\left\lbrack\mathbf{1}_A\vert \Omega(s)^+\right\rbrack - P_{p}(\Omega(s)^-)E_{p}\left\lbrack\mathbf{1}_A\vert \Omega(s)^-\right\rbrack \\
&=& (p' - p) E_p\left\lbrack \mathbf{1}_A - \mathbf{1}_A\circ j_s\vert \omega\in\Omega(s)^+\right\rbrack
+ E_{p',s}\left\lbrack\mathbf{1}_A\vert \Omega(s)^-\right\rbrack - E_{p}\left\lbrack\mathbf{1}_A\vert \Omega(s)^-\right\rbrack \\
&=& (p' - p) E_p\left\lbrack  \mathbf{1}_{\mathrm{Piv}(A,s)\cap\Omega(s)^+}\vert \omega\in\Omega(s)^+\right\rbrack \\
&=& (p' - p) E_p\left\lbrack  \mathbf{1}_{\mathrm{Piv}(A,s)}\right\rbrack = (p' - p) P_p(\mathrm{Piv}(A,s)) .
\end{eqnarray*}
That implies finally
$$
\frac{\partial P_p (A)}{\partial p(s)} = P_p(\mathrm{Piv}(A,s)).
$$
{\rm (iv)} By $A\in\FF_{\GG}$, we have
$$
E_{p}\mathbf{1}_A = E_{p}\left\lbrack \mathbf{1}_A\vert \FF_{\GG}\right\rbrack = \sum_{\omega\in\lbrace 0,1\rbrace^{\GG}}\Pi_{s\in \GG} 1_A(\omega)\Pi_{s\in \GG}P_{p}(\omega(s))
$$
that means, we may think of the distribution $P_p$ as a distribution depending on {\em finitely many} real parameters $\lbrace p(s)\,:\,s\in\GG\rbrace$. That implies together with {\rm (iii)}
$$
\frac{d}{dp}P_p(A) = \sum_{s\in \GG}\frac{\partial P_p (A)}{\partial p(s)}\frac{\partial p(s)}{\partial p} = \sum_{s\in \GG}P_p(\mathrm{Piv}(A,s)) = \sum_{s\in\GG} E_p\mathbf{1}_{\mathrm{Piv}(A,s)} = E_pN_A .
$$
\end{Proof}

\subsection{The proof of (\ref{EXPECTBOUND})}\label{pf}

We follow closely the proof in \cite{Gri:99}, p. 263 ff. Let $P_p(x,y)$ denote the probability that there is a path connecting the sites $x$ and $y$ and $P_p^{(N)}(x,y)$ the probability that there is a path connecting $x$ and $y$ which lies entirely in $\TT^{(N)}$. Now
$$
\chi_N(p,y) := \sum_{x\in\TT^{(N)}} P_p^{(N)}(x,y)
$$
is the expected size of the connected cluster around $y$ in $\TT^{(N)}$ and
$$
\chi (p,y) := \sum_{x\in\TT} P_p(x,y)
$$
the expected cluster size in $\TT$. Note that $\chi (p) = \chi(p,0)$. Furthermore, we write
$$
\chi_N(p) := \max \lbrace \chi_N(p,y)\,:\,y\in\TT^{(N)}\rbrace .
$$
\noindent {\rm (i)} First of all,
$$
\chi (p) \geq \chi_N(p) \geq\chi_N(p,0) =\sum_{x\in\TT^{(N)}} P_p^{(N)}(x,0) \to \sum_{x\in\TT} P_p(x,y) = \chi(p)
$$
implies that we have by bounded convergence
$$
\lim_{N\to\infty} \chi_N(p) = \chi(p).
$$
{\rm (ii)} Denote by $A_N(x,y)$ the event that there is a path connecting $x$ and $y$ in $\TT^{(N)}$. Then, by Russo's formula,
$$
\frac{d}{dp} \chi_N(p,y) = \sum_{x\in\TT^{(N)}} \sum_{s\in\TT^{(N)}} P_p(\mathrm{Piv}(A_N(x,y),s)).
$$
A site $s\in\TT^{(N)}$ is now pivotal for $A_N(x,y)$, if and only if
\begin{enumerate}
\item $s$ is adjacent to two different and non - adjacent sites $x'$ and $y'$.
\item There is a path connecting $x$ and $x'$.
\item There is a {\em disjoint} path connecting $y$ and $y'$, meaning that no site in this path is adjacent to any site in the path connecting $x$ ans $x'$.
\end{enumerate}
This means that switching $s$ on or off will switch a connection between $x$ and $y$ on or off (which changes the value of the corresponding indicator function). Having disjoint paths between different pairs of sites is a typical example of a disjointly occuring event and therefore we can write the three conditions above shortly by saying that for all $x,y\neq z\in\TT^{(N)}$ and all $x'\neq y'$ adjacent to and different from $s$, we have
$$
A_N(x,x')\circ A_N(y,y')\subset\mathrm{Piv}(A_N(x,y),s)
$$
and that on the other hand
$$
\mathrm{Piv}(A_N(x,y),s) = \bigcup_{x'\neq y'\,\mathrm{adjacent \,\,to}\,\, s} A_N(x,x')\circ A_N(y,y').
$$
That implies by BK inequality
$$
P_p(\mathrm{Piv}(A_N(x,y),s)) \leq \sum_{x'\neq y'\,\mathrm{adjacent \,\,to}\,\, s} P_p^{(N)}(x,x') \,P_p^{(N)}(y,y').
$$
Finally inserting this into Russo's formula yields
\begin{eqnarray*}
\frac{d}{dp} \chi_N(p,y) &=& \sum_{x\in\TT^{(N)}} \sum_{s\in\TT^{(N)}} P_p(\mathrm{Piv}(A_N(x,y),s)) \\
&\leq& \sum_{x\in\TT^{(N)}} \sum_{s\in\TT^{(N)}} \,\,\,\sum_{x'\neq y'\,\mathrm{adjacent \,\,to}\,\, s} P_p^{(N)}(x,x') \,P_p^{(N)}(y,y') \\
&=&  \sum_{s\in\TT^{(N)}} \,\,\,\sum_{x'\neq y'\,\mathrm{adjacent \,\,to}\,\, s} \chi_N(p,x') \,P_p^{(N)}(y,y') \\
&\leq&  \chi_N(p)\sum_{s\in\TT^{(N)}} \,\,\,\sum_{x'\neq y'\,\mathrm{adjacent \,\,to}\,\, s} P_p^{(N)}(y,y') \\
&=&  3\,\chi_N(p)\sum_{s\in\TT^{(N)}} \,\,\,\sum_{y'\,\mathrm{adjacent \,\,to}\,\, s} P_p^{(N)}(y,y') \\
&=&  3\times 6\,\chi_N(p)\sum_{s\in\TT^{(N)}}  P_p^{(N)}(y,s) \\
&=&  18\,\chi_N(p)\,\chi_N(p,y) \leq 18 \,\chi_N(p)^2 .
\end{eqnarray*}
{\rm (iii)} Integrating this differential inequality over the interval $\lbrack p, p_c\rbrack$ yields
$$
\frac{1}{\chi_N(p)} - \frac{1}{\chi_N(p_c)} \leq 18\,(p-p_c)
$$
(for details, see the above mentioned proof in \cite{Gri:99}) and by {\rm (i)} $\chi_N \to \chi$ and the fact that $\chi (p_c) = \infty$ we finally obtain
$$
\chi (p) \geq \frac{1}{18 \,(p-p_c)}.
$$

\subsection{Matching graphs and {$p_c = 1/2$}}

\noindent In this subsection, we will shortly review the material from \cite{SykEss:64} about site percolation and matching graphs\index{matching graphs}. We start with a finite graph $\GG$ with $N$ sites. The probability that a site is marked {\em active} (or {\em black}) is given by $p$, the probability that it is marked {\em inactive} (or {\em white}) is $q = 1- p$. Denote a connected cluster of black points by $C$ and its boundary by
$$
\partial C := \lbrace s\in \GG - C\,:\,s \,\,\mathrm{is\,\, adjacent\,\, to \,\, some \,\,site} \,\, s'\in C\rbrace .
$$
That means, the expected cluster size is a polynomial in $p$ and $q$ given by
$$
K(p,q,\GG) = E\,\vert C \vert = \sum_{C\subset \GG} \vert C\vert\, p^{\vert C\vert}q^{\vert \partial C\vert}.
$$
By reversing the roles of $p$ and $q$, we obtain the expected numbers of white clusters. To extend this concept to infinite graphs, we consider the {\em mean cluster size per site}
\begin{equation}\label{poly}
k(p,q,\GG) = E (\vert C \vert / \vert\GG\vert ) = \frac{1}{\vert \GG\vert}\sum_{C\subset \GG} \vert C\vert\, p^{\vert C\vert}q^{\vert \partial C\vert},
\end{equation}
use a proper exhaustion $\GG_1\subset\GG_2\subset ... \subset\GG$ of an infinite graph $\GG$ and consider the formal power series
$$
k(p,q,\GG) = \lim_{k\to\infty} \frac{1}{\vert \GG_k \vert}\sum_{C\subset \GG_k} \vert C\vert\, p^{\vert C\vert}q^{\vert \partial C\vert}
$$
which shows that we obtain in this case the expected {\em finite} cluster size per size, taking into account only finite subclusters from $\GG$. By
\begin{equation}\label{lowhi}
\begin{array}{ll}  k_L(p,\GG) = k (p, 1-p,\GG), &  k_H(q,\GG) = k (1 - q, q,\GG),\end{array}
\end{equation}
we clearly have $k_L(p,\GG)=k_H(1-p,\GG)$ and $k_H(q,\GG)=k_L(1-q,\GG)$.

\begin{Definition} We call a (possibly infinite) graph $\GG$ {\em self - matching}, if there is a polynomial $\varphi_{\GG}(p)$ such that
\begin{equation}\label{matching}
k_L(p,\GG) = \varphi_{\GG}(p) + k_H(p,\GG) .
\end{equation}
$\varphi_{\GG}$ is called the {\em matching polynomial}.
\end{Definition}

\begin{Theorem} The triangular lattice is self-matching with
$$
\varphi_{\TT} (p) = p - 3p^2 + 2 p^3.
$$
\end{Theorem}

\begin{Proof} See \cite{SykEss:64}.\end{Proof}

\noindent When we assume that $k_L(p,\TT)$ has precisely one pole at the {\em critical percolation probability} $p_c$ for the triangular lattice (see for instance \cite{Kes:82}), we can actually use the preceding theorem to determine $p_c$. Here, the special form of the matching polynomial does not play any role, only the fact that it is a polynomial and hence bounded on $p\in\lbrack 0,1\rbrack$ is important. Therefore $k_L(p_c,\TT)=\infty$ implies $k_H(p_c,\TT)= k_L(1 - p_c,\TT) = \infty$. The assumption that there is only one pole immediately implies $p_c = 1- p_c$ and thus $p_c=1/2$.\\

\begin{Remark} If the graph $\GG$ is not {\em self - matching}\index{self-matching}, we can construct a so called {\em matching graph} $\GG^*$ (for the construction, see again \cite{SykEss:64}, or \cite{Kes:82}) with the same number of vertices such that  instead of (\ref{matching}), we have
$$
\begin{array}{ll}  k_L(p,\GG) = \varphi (p) + k_H(p,\GG^*), &   k_L(p,\GG^*) = \varphi^* (p) + k_H(p,\GG),  \end{array}
$$
together with some relations between $\varphi$ and $\varphi^*$ and these equations can be used in a similar way as above to obtain some information about the critical probability. $(\GG,\GG^*)$ is called a {\em matching pair}. For self - matching graphs, we have $\GG^* = \GG$.
\end{Remark}

\section[Bounded detector devices]{Bounded detector devices}

In the discussion of realistic signals, we introduced the notion of a {\em bounded detector device}. In statistical terminology, those devices from an instance of the method of truncation. A {\em bounded detector device of range }\index{bounded detector device} $r>0$ is only able to display signal strengths $Y(s)$ with intensities between $-r$ and $r$. Thus, the effect of the detector device on a signal $Y$ is that instead of the full information about $Y(s)$, $s\in S$, only the information contained in the cutoff signal
\begin{equation}\label{bdd}
D(Y) = {\mathrm{max}} \lbrace {\mathrm{min}} \lbrace Y, r \rbrace ,-r \rbrace
\end{equation}
is used for further analysis. Intensities of absolute value larger than $r$ can simply not be registered and all information about the behavior of the signal above and below the cutoff is lost before the signal processing even starts.\\

The detection results in the present paper were proved for bounded signals. What happens if this assumptions doesn't hold? First of all, from a purely mathematical point of view, the notion of bounded detector devices can be equivalently reformulated by saying that all considerations are only valid as statements that are obtained while {\em conditioning} on the event $\lbrace \vert Y\vert < r\rbrace$. In other words, all results are still valid without any change if we understand them as being obtained by {\em conditioning on the event}
\begin{equation}\label{conditioning}
D_0 := \lbrace D(Y)=Y \rbrace.
\end{equation}
Of course, the probability $\pi_D := P(D_0)$ then yields an important characteristic of the detector device, and it could be often desirable to have $\pi_D$ close to one. However, a deeper analysis of biological, engineering and cybernetical aspects of the problem leads us to the following extremely useful observation.




We think of signal processing as consisting of at least three different parts,
\begin{enumerate}
\item a {\em filter} which has the purpose to transform the incoming signal to fit in an optimal way into the bounds of the detector device,
\item the {\em bounded detector device} as described above, and
\item the {\em processor}, which analyses the detected signal $D(Y)$ and determines what is finally perceived.
\end{enumerate}
We thus arrive at the following scheme
$$
\begin{array}{lllll}\mathrm{Signal} \rightarrow  &\fbox{Filter}\rightarrow&\fbox{Detector}\rightarrow&\fbox{Processor}\rightarrow&\mathrm{Perception}\end{array}
$$
where the detector is the fixed component, the filter is chosen on the basis of the incoming signal and the bounds of the detector and the processor algorithm is chosen on the basis of knowledge about the detector and the chosen filter. Choosing an appropriate filter for a given environment is thus another problem of perception, a problem that we will not address in these notes.\\

\noindent{\bf Example.} As an example, as the {\em detector device} of the human eye, we only consider the photo receptors situated at the retina, the processor obviously is the brain, and the filter is given by lens and iris which adapt to different light intensities for instance in night vision, but can also be those parts together with another device like, for instance, sun glasses.\\

For a visual perception of any system in biology or cybernetics, the meaning of a \emph{good} Filter is exactly to filter out (or to transform) the incoming information in such a way that the Detector might still perceive a reasonable part of reality, despite the fact that the Detector works with signals in the diapason $[-r, \, r]$ only. Say, in the above Example, a human eye doesn't have to properly perceive ultraviolet and infrared light frequencies in order to be able to see trees. A human brain doesn't need to process any information that could come with ultraviolet and infrared lights either.

This implies that our consideration of bounded detector devices fits many important biological situations. Moreover, working with bounded detector devices can be profitable for construction of artificial vision systems in robotics. A robot needs to perceive and to process only signals and information within the diapason that fits his tasks.

\end{document}